\pdfoutput=1
\RequirePackage{ifpdf}
\ifpdf 
\documentclass[pdftex]{sigma}
\else
\documentclass{sigma}
\fi

\numberwithin{equation}{section}

\newtheorem*{Theorem*}{Theorem}

 { \theoremstyle{definition}

 }

\newcommand\mytop[2]{\genfrac{}{}{0pt}{}{#1}{#2}}
\def\genhyperF#1#2#3#4#5#6{{{}_{#1}F_{#2}}\!\left(\genfrac{}{}{0pt}{}{#3, #4}{#5};#6\right)}
\def\expe{{\rm e}}

\def\iunit{{\rm i}}
\def\id{\,{\rm d}}
\def\bigO{{\mathcal O}}

\begin{document}
\allowdisplaybreaks

\renewcommand{\thefootnote}{}

\newcommand{\arXivNumber}{2311.11886}

\renewcommand{\PaperNumber}{023}

\FirstPageHeading

\ShortArticleName{Lerch $\Phi$ Asymptotics}

\ArticleName{Lerch $\boldsymbol{\Phi}$ Asymptotics\footnote{This paper is a~contribution to the Special Issue on Asymptotics and Applications of Special Functions in Memory of Richard Paris. The~full collection is available at \href{https://www.emis.de/journals/SIGMA/Paris.html}{https://www.emis.de/journals/SIGMA/Paris.html}}}

\Author{Adri B.~OLDE DAALHUIS}

\AuthorNameForHeading{A.B.~Olde Daalhuis}

\Address{School of Mathematics and Maxwell Institute for Mathematical Sciences,\\ The University of Edinburgh, Edinburgh EH9~3FD, UK}
\Email{\href{mailto:A.OldeDaalhuis@ed.ac.uk}{A.OldeDaalhuis@ed.ac.uk}}
\URLaddress{\url{http://www.maths.ed.ac.uk/~adri/}}

\ArticleDates{Received November 22, 2023, in final form March 11, 2024; Published online March 21, 2024}

\Abstract{We use a Mellin--Barnes integral representation for the Lerch transcendent $\Phi(z,s,a)$ to obtain large $z$ asymptotic approximations. The simplest divergent asymptotic approximation terminates in the case that $s$ is an integer. For non-integer $s$ the asymptotic approximations consists of the sum of two series. The first one is in powers of~$(\ln z)^{-1}$ and the second one is in powers of $z^{-1}$. Although the second series converges, it is completely hidden in the divergent tail of the first series. We use resummation and optimal truncation to make the second series visible.}

\Keywords{Hurwitz--Lerch zeta function; analytic continuation; asymptotic expansions}

\Classification{11M35; 30E15; 41A30; 41A60}

\renewcommand{\thefootnote}{\arabic{footnote}}
\setcounter{footnote}{0}

\section{Introduction and summary}
This paper is in memory of Richard Paris. He was an expert in special functions and
asymptotics, especially the use of Mellin--Barnes integral representations.
Most of the techniques that we use in this paper are all coming from his
book \cite{Paris2001} (co-authored by David Kaminski). Note that in~\cite{paris2016stokes} Richard
discusses the large $a$ asymptotics for the Lerch transcendent, and he uses a~Mellin--Barnes
integral representation that is with respect to $a$, whereas \eqref{MBint} is with respect to~$z$.

The Lerch transcendent is defined via
\begin{equation}\label{Phi}
 \Phi(z,s,a)=\sum_{n=0}^{\infty}\frac{z^{n}}{(a+n)^{s}},
 \qquad |z|<1,
\end{equation}
and via analytic continuation (see \eqref{int}) elsewhere in the complex $z$ plane.
If $s$ is not an integer then $|\arg a|<\pi$; if $s$ is a positive integer then
$a\not= 0,-1,-2,\ldots$; if $s$ is a non-positive integer then $a$ can be any complex number.
The Lerch transcendent has a branch-point at $z=1$, which follows from the expansion
(see \cite[Section~1.11\,(8)]{Erdelyi:1953:HTF1})
\begin{equation*}
 \Phi(z,s,a)=\Gamma(1-s) z^{-a} (-\ln z)^{s-1}+z^{-a}
 \sum_{n=0}^{\infty}\frac{\zeta(s-n,a)}{n!}(\ln z)^n,
\end{equation*}
$|\ln z|<2\pi$, $s\not=1,2,3,\ldots$, $a\not=0,-1,-2,\ldots$, in which
$\zeta(s,a)=\Phi(1,s,a)$ is the Hurwitz zeta function.
The principle branch for $\Phi(z,s,a)$ is the sector $|\arg (1-z)|<\pi$.

In several publications the main tool to obtain asymptotic expansions is the integral
representation
\begin{equation}\label{int}
 \Phi(z,s,a)=\frac1{\Gamma(s)}\int_0^\infty
 \frac{x^{s-1}\expe^{-ax}}{1-z\expe^{-x}}\id x,\qquad \operatorname{Re} s>0,\quad\operatorname{Re} a>0.
\end{equation}
The small and large $a$ asymptotics is discussed in \cite{Cai2019,Ferreira2004,paris2016stokes} (and \cite{Nemes2017} for the case $z=1$), and
the asymptotics as $\operatorname{Re} s\to -\infty$ is discussed in \cite{Navas2013}.
In \cite{Ferreira2004}, the large $z$ asymptotics is also discussed, and the authors obtain
\begin{equation}\label{FL1}
 \Phi(-z,s,a)\sim\frac1{\Gamma(s)}\sum_{n=1}^\infty
 \frac{A_n(z,s,a)}{(-z)^{n}}+\frac{(\ln z)^s}{z^a\Gamma(s)}
 \sum_{n=0}^\infty\frac{B_n(s,a)}{(\ln z)^{n+1}}
\end{equation}
as $z\to\infty$, with
\begin{equation}\label{FL2}
 A_n(z,s,a)=\frac{\Gamma(s,(a-n)\ln z)-\Gamma(s)}{(a-n)^s}
 =\frac{-(\ln z)^s}{s}M(s,s+1;(a-n)\ln z),
\end{equation}
and
\begin{gather*}
 B_0(s,a)=\frac{\psi\bigl(\frac{a+1}2\bigr)-
 \psi\bigl(\frac{a}2\bigr)}2,\\
 B_n(s,a)=\frac{n!{s-1\choose n}}{2^{n+1}}\left( \zeta\left(n+1,\frac{a}2\right)-
 \zeta\left(n+1,\frac{a+1}2\right)\right),
\end{gather*}
$n=1,2,3,\ldots$. The second representation for $A_n(z,s,a)$ in \eqref{FL2} is in terms
of the Kummer $M$ confluent hypergeometric function
(see \cite[\href{https://dlmf.nist.gov/13.2.E2}{equation~(13.2.2)}]{NIST:DLMF}),
and is more convenient from an analytical/numerical point of view.
In \cite{Ferreira2004}, they do truncate the two sums in \eqref{FL1} and give sharp error bounds
for each truncated sum. It seems not possible to combine the two error bounds
to obtain an approximation
that is correct up to a combined error estimate of, say, $\bigO\bigl(z^{-N}\bigr)$,
especially because the second sum in \eqref{FL1} is divergent, and its terms completely
dominate the combined asymptotics. The incomplete gamma function in the first series seems
misplaced, because for each $n$ we have
\smash{$\frac{\Gamma(s,(a-n)\ln z)}{(a-n)^s (-z)^{n}}
=\bigO\bigl(\frac{(\ln z)^s}{z^a}\bigr)$} as $z\to\infty$.
Hence, it seems that the incomplete gamma function should be incorporated in the second
series of \eqref{FL1}. However, we do observe that the advantage of definition
\eqref{FL2} is that in this way the $A_n(z,s,a)$ are bounded as $a\to n$.

Below, we will show that the large $z$ asymptotics simplifies dramatically
when $s$ is an integer. It will take a considerable amount of work to
obtain that result from \eqref{FL1}.

It is very surprising that it is not well known that the Lerch transcendent has a simple
Mellin--Barnes integral representation. This is a more obvious tool to obtain asymptotics.
For the moment we take $\operatorname{Re} a>0$. In this case,
we have the Mellin--Barnes integral representation (see \cite[equation~(3.11)]{Srivastava2011})
\begin{equation}\label{MBint}
 \Phi(-z,s,a)=\frac1{2\pi\iunit}\int_L
 \frac{\Gamma(1+t)\Gamma(-t)z^t}{(a+t)^s}\id t,
\end{equation}
in which $L$ is a contour from $-\iunit\infty$ to $\iunit\infty$ which crosses the
real $t$ axis somewhere in the interval~$(\max(-\operatorname{Re} a,-1),0)$. When we push this contour to the
right we will get contributions from the poles of $\Gamma(-t)$ and this results in
expansion~\eqref{Phi}. Integral representation~\eqref{MBint} supplies an
analytic continuation to the sector $|\arg z|<\pi$.

We write for the moment $z=\rho\expe^{\iunit\theta}$ and observe that for the dominant factor
of the integrand in~\eqref{MBint} we have
$\bigl|\Gamma(1+t)\Gamma(-t)z^t\bigr|\sim 2\pi \rho^{\operatorname{Re} t}\expe^{(\mp\pi-\theta)\operatorname{Im} t}$
as $\operatorname{Im} t\to\pm\infty$. The factor $\rho^{\operatorname{Re} t}$ can be used to analytically continue
integral representation~\eqref{MBint}.
In the case $0<\rho<1$, we bent contour $L$ to the right (say we go from $(1-\iunit)\infty$
to $(1+\iunit)\infty$) and with this choice for~$L$ integral representation~\eqref{MBint} is also valid for $z\in(-1,0)$.
In the case that $\rho>1$, we bent contour~$L$ to the
left (say we go from $(-1-\iunit)\infty$ to $(-1+\iunit)\infty$) and
with this choice for $L$ integral representation~\eqref{MBint}
is valid across the $z$-branch-cut $(-\infty,-1)$.

To obtain a large $z$ asymptotic expansion all that we have to do is to push $L$ to the
left. This is what they do in \cite{Srivastava2011}, but somehow they do miss the fact
that there will also be contributions from the branch-point at $t=-a$.

In Section~\ref{Sect:largez}, we will study the large $z$ asymptotics by pushing contour $L$ to the left.
The contributions of the poles of $\Gamma(1+t)$ will give us a simple convergent asymptotic
expansion in powers of $z^{-1}$ that can be expressed in terms of a Lerch transcendent.
The contribution of the branch-point at $t=-a$ is more complicated. In the case that $s=S$ is an integer,
this contribution is just a finite sum in powers of $1/\ln z$, and we obtain
\begin{gather}\label{BIsimple}
 \Phi(-z,s,a)=\frac{2\pi\iunit}{z^{a}}\sum_{n=0}^{S-1}
 \frac{b_n}{\Gamma(S-n)(\ln z)^{n-S+1}}
 -\sum_{n=1}^\infty\frac{(-z)^{-n}}{\left(a-n\right)^S},
\end{gather}
in which the first sum is zero in the case that $S$ is a non-positive integer. The coefficients
$b_n$ are the Taylor coefficients of $\frac{\iunit}2 \csc \pi(t-a)$ about $t=0$.
This result has already been used in \cite[Lemma 3]{NCA2023}.

In the case that $s$ is not an integer, the first series in \eqref{BIsimple} does not terminate
and it is a~divergent series. We will use resummation and optimal truncation to obtain
an approximation for $\Phi(-z,s,a)$ that is accurate up to order
$\bigO\left(z^{-{N-1}}\right)$:

\begin{theorem}\label{thm1}
 We take $|\arg(1+z)|\leq \pi$, $a$ and $s$ bounded complex numbers,
 $\operatorname{Re} a>0$. Let~$N$ be a~fixed integer with $N>\operatorname{Re} a$. Take $|z|$ large and let $M$
 be a positive integer such that $|M-s+1|\approx |(N+1-a)\ln z|$. Then
 \begin{align*}
 \Phi(-z,s,a)={}& \sum_{n=0}^N \frac{(-z)^n \Gamma(s,(a+n)\ln z)}{(a+n)^s\Gamma(s)}
 +\frac{2\pi\iunit}{z^{a}}\sum_{m=0}^{M-1}
 \frac{b_{m,N}}{\Gamma(s-m)(\ln z)^{m-s+1}}\\
 & +\frac1{\Gamma(s)}\sum_{n=1}^N \frac{A_n(z,s,a)}{%
 (-z)^n}+\bigO\bigl(z^{-N-1}\bigr),
\end{align*}
as $z\to\infty$ uniformly with respect to $\arg z\in[-\pi,\pi]$.
The $A_n(z,s,a)$ are defined in \eqref{FL2} and
\begin{align}
 b_{0,N}={}&\frac{(-1)^N}{4\pi\iunit}\left(\psi\left(\frac{N+1+a}2\right)
 -\psi\left(\frac{N+2+a}2\right)\right.\nonumber\\
 &\left.-\psi\left(\frac{N+1-a}2\right)
 +\psi\left(\frac{N+2-a}2\right)\right),\nonumber\\
 b_{n,N}={}&\frac{(-1)^{N}}{2^{n+2}\pi\iunit}\left(\zeta\left(n+1,\frac{N+2+a}2\right)
 -\zeta\left(n+1,\frac{N+1+a}2\right)\right.\nonumber\\
 &\left.+(-1)^n\zeta\left(n+1,\frac{N+1-a}2\right)
 -(-1)^n\zeta\left(n+1,\frac{N+2-a}2\right)\right),\label{BnN2}
\end{align}
$n=1,2,3,\ldots$.
\end{theorem}

 We numerically verify this result in Table \ref{tab:combined}.
Observe above that the $M$ does depend on~$z$. In the $m$-series we do take an optimal
number of terms. Hence, its remainder will be `exponentially-small'.
In the proof below we do show that it is $\bigO\bigl(\expe^{-(N+1)|x|}\bigr)$,
with $x=\ln z$.

Often when one encounters divergent Poincar\'e asymptotic series, it is possible to convert it
to a convergent factorial series. In the case that $s$ is not an integer it would be convenient
to replace the first series in~\eqref{BIsimple} (which diverges) by a convergent factorial
series and in that way we can also incorporate in our approximation the full second series
of~\eqref{BIsimple}. We do create a~convergent factorial-type series in
Section~\ref{Sect:factorial}. As usual with these types of series, the convergence is very slow, and
it is not easy to obtain sharp error estimates.

\section[Large z asymptotics]{Large $\boldsymbol{z}$ asymptotics}\label{Sect:largez}
We start with $0<a<1$ and we can use analytic continuation afterwards.
We take $|z|>1$ in the sector $|\arg(1+z)|\leq\pi$ and push $L$ to the left and obtain
\begin{gather}\label{BI}
 \Phi(-z,s,a)=B(z,s,a)+I(z,s,a).
\end{gather}
The first term is the contribution of the branch-point at $t=-a$
\begin{gather}\label{Bcontr}
 B(z,s,a)=\frac1{2\pi\iunit}\int_{-(1+\iunit)\infty}^{(-a+)}
 \frac{\Gamma(1+t)\Gamma(-t)z^t}{(a+t)^s}\id t,
\end{gather}
in which the contour begins at $-(1+\iunit)\infty$, circles $t=-a$ once in the positive direction,
and returns to $-(1+\iunit)\infty$. We did observe in the second paragraph below \eqref{MBint}
that in the case~$|z|>1$ the integrals converge across the $z$-branch cuts $\arg z=\pm\pi$.
The second term in \eqref{BI} is the sum of the residue contributions of the poles of $\Gamma(1+t)$
\begin{gather}\label{Icontr}
 I(z,s,a)=-\expe^{-\pi\iunit s}\sum_{n=1}^\infty\frac{(-z)^{-n}}{(n-a)^s}
 =\expe^{-2\pi\iunit s}a^{-s}-\expe^{-\pi\iunit s}\Phi (-1/z,s,-a ).
\end{gather}
The infinite series representation of $I(z,s,a)$ is already a large $z$ asymptotic
expansion. Hence, all we need is an asymptotic expansion for $B(z,s,a)$. We
modify integral representation~\eqref{Bcontr}~as%
\begin{gather}\label{Bcontr2}
 B(z,s,a)=z^{-a}\int_{-(1+\iunit)\infty}^{(0+)}
 \expe^{t\ln z} t^{-s}g(t)\id t,
\end{gather}
with
\begin{gather}\label{g}
 g(t)=\frac{\iunit}{2\sin\pi(t-a)}
 =\lim_{N\to\infty} \frac{\iunit}{2\pi}\sum_{n=-N}^N\frac{(-1)^n}{t-a-n}.
\end{gather}
We use this sum representation of $g(t)$ in \eqref{Bcontr2} and obtain the expansion
\begin{gather}\label{Bgamma}
 B(z,s,a)=\lim_{N\to\infty} \sum_{n=-N}^N \frac{(-z)^n}{(a+n)^s\Gamma(s)}
 \Gamma(s,(a+n)\ln z),
\end{gather}
in terms of the incomplete gamma function. This expansions converges slowly, but has no asymptotic
property as $z\to\infty$.

To obtain a simple Poincar\'e asymptotic expansions in inverse powers of $\ln z$,
we have to expand $g(t)$ about the origin. Let
\begin{gather}\label{gMacl}
 g(t)=\sum_{n=0}^\infty b_nt^n,\qquad b_0=\frac{1}{2\iunit\sin\pi a},\qquad
 b_1=2\pi\iunit b_0^2\cos\pi a.
\end{gather}
The reader can check that $g(t)g''(t)=2g'^2(t)+\pi^2g^2(t)$, from which we obtain the recurrence relation{\samepage
\begin{align*}
 (n+2)(n+1)b_0b_{n+2}={}&\sum_{m=0}^n\big(2(m+1)(n-m+1)b_{m+1}b_{n-m+1}
 +\pi^2b_m b_{n-m}\big)\\
 &-\sum_{m=0}^{n-1}(m+2)(m+1)b_{m+2}b_{n-m},\qquad n\geq0.
\end{align*}
 Hence, the computation of the coefficients is straightforward.}

To obtain an asymptotic expansion for the right-hand side of \eqref{Bcontr2},
we will use Watson's lemma for loop integrals. See \cite[Section~4.5.3]{Olver1997}.
We substitute \eqref{gMacl} into \eqref{Bcontr2} and obtain the asymptotic expansion
\begin{equation}\label{Basymp}
 B(z,s,a)\sim \frac{2\pi\iunit}{z^{a}}\sum_{n=0}^\infty
 \frac{b_n}{\Gamma(s-n)(\ln z)^{n-s+1}},
\end{equation}
as $z\to\infty$.

In the case that $s=S$ is an integer, the infinite series on the right-hand side
of \eqref{Basymp} terminates.
If $S$ is a non-positive integer, then $B(z,S,a)=0$, and if $S$ is a positive integer,
we have
\begin{gather*}
 B(z,S,a)= \frac{2\pi\iunit}{z^{a}}\sum_{n=0}^{S-1}
 \frac{b_n}{\Gamma(S-n)(\ln z)^{n-S+1}}.
\end{gather*}

However, in the case that $s$ is not an integer
\eqref{Basymp} is a divergent asymptotic expansion. For modestly large $|z|$, the $\ln z$ in
this expansion is not very large, and the optimal number of terms will be small.
In Section~\ref{Sect:factorial}, we will obtain another type of asymptotic approximation for~$B(z,s,a)$.
Here we are going to combine the three expansions \eqref{Basymp}, \eqref{Icontr} and \eqref{Bgamma}.

First, we observe that when we take an optimal number of terms in \eqref{Basymp} the smallest term will
still be bigger than the first term in \eqref{Icontr}. The optimal number of terms is connected to the
distance between the origin and the nearest pole of $g(t)$. By incorporating in the expansion
the poles of $g(t)$ that are near the origin, we `slow down' the divergence of \eqref{Basymp}.

\begin{proof}[Proof of Theorem \ref{thm1}]
In \eqref{BI}, we write $\Phi(-z,s,a)=B(z,s,a)+I(z,s,a)$ and for $I(z,s,a)$
we have the convergent asymptotic expansion \eqref{Icontr}. Hence, we need an asymptotic
approximation for $B(z,s,a)$ with the correct remainder estimate.
Let $N$ be a fixed positive integer and take
\begin{gather*}
 g_N(t)=g(t)-\frac{\iunit}{2\pi}\sum_{n=-N}^N\frac{(-1)^n}{t-a-n},
\end{gather*}
compare \eqref{g}. Then
\begin{gather*}
 B(z,s,a)=\sum_{n=-N}^N \frac{(-z)^n}{(a+n)^s\Gamma(s)}\Gamma(s,(a+n)\ln z)
 +z^{-a}\int_{-(1+\iunit)\infty}^{(0+)}\expe^{t\ln z} t^{-s}g_N(t)\id t.
\end{gather*}
We have
\begin{gather*}
 g_N(t)=\sum_{n=0}^\infty b_{n,N}t^n,
\end{gather*}
with
\begin{gather}\label{BnN}
 b_{n,N}=b_n+\frac{\iunit}{2\pi}\sum_{m=-N}^N \frac{(-1)^m}{(a+m)^{n+1}}.
\end{gather}
A numerical more stable presentation for the $b_{n,N}$ is \eqref{BnN2}.

For large $n$, the $b_{n,N}$ will be of the size of the first omitted terms in \eqref{BnN},
that is, for the case~$0<a<1$, the term with $m=-N-1$.
Thus
\smash{$b_{n,N}\sim \frac{\iunit}{2\pi}(-1)^N (a-N-1)^{-n-1}$} as $n\to\infty$.
We will estimate the remainder in
\begin{gather*}
 z^{-a}\int_{-(1+\iunit)\infty}^{(0+)}\expe^{t\ln z} t^{-s}g_N(t)\id t
 =\frac{2\pi\iunit}{z^{a}}\sum_{m=0}^{M-1}
 \frac{b_{m,N}}{\Gamma(s-m)(\ln z)^{m-s+1}}+R_M^{(B)}(z)
\end{gather*}
via the machinery in \cite[Section~4]{OD98}, because the integrals are exactly of the same form
as the ones discussed in that paper.
We want the remainder after taking the optimal number of terms. Hence, $M$ will be large.
We have
\begin{align}
 \left|R_M^{(B)}(z)\right|&=\bigO\left(\frac{z^{-a} (\ln z)^{s-M}}{%
 (N+1-a)^{M}}F^{(1)}\left((a-N-1)\ln z;\mytop{M-s+1}{1}\right)\right)\nonumber\\
 &=\bigO\left(\frac{z^{-a} \Gamma(M-s+1)\sqrt{M-s+1}}{%
 \left((N+1-a)\ln z\right)^{M-s+1}}\right),\label{RM1}
\end{align}
in which both $M$ and $z$ are large. For the first hyperterminant $F^{(1)}$ see
\cite[Appendix A]{BHNOD2018}. The second equal sign in \eqref{RM1} follows from
\cite[Proposition~B.1]{BHNOD2018}.
It follows that for the optimal~${M=M_{\rm opt}}$, we have
$|M_{\rm opt}-s+1|\sim |(N+1-a)\ln z|$. With this choice for $M$, we~have
\begin{gather*}
 \bigl|R_M^{(B)}(z)\bigr|=\bigO\left(\frac{z^{-a}
 \Gamma(M-s+1)}{(M-s+1)^{M-s+\frac12}}\right)=
 \bigO\bigl(z^{-a}\expe^{s-M-1}\bigr)=\bigO\bigl(z^{-N-1}\bigr),
\end{gather*}
as $z\to\infty$.

Combining the results above, we have the approximation
\begin{align}
 \Phi(-z,s,a)={}&
 \sum_{n=-N}^N \frac{(-z)^n \Gamma(s,(a+n)\ln z)}{(a+n)^s\Gamma(s)}
 +\frac{2\pi\iunit}{z^{a}}\sum_{m=0}^{M_{\rm opt}-1}
 \frac{b_{m,N}}{\Gamma(s-m)(\ln z)^{m-s+1}}\nonumber\\
 & -\expe^{-\pi\iunit s}\sum_{n=1}^N\frac{(-z)^{-n}}{(n-a)^s}
 +R_N(z),\nonumber\\
 ={}& \sum_{n=0}^N \frac{(-z)^n \Gamma(s,(a+n)\ln z)}{(a+n)^s\Gamma(s)}
 +\frac{2\pi\iunit}{z^{a}}\sum_{m=0}^{M_{\rm opt}-1}
 \frac{b_{m,N}}{\Gamma(s-m)(\ln z)^{m-s+1}}\nonumber\\
 & +\frac1{\Gamma(s)}\sum_{n=1}^N \frac{A_n(z,s,a)}{(-z)^n}+R_N(z),\label{Basympcombined}
\end{align}
with $R_N(z)=\bigO\left(z^{-N-1}\right)$, as $z\to\infty$ uniformly with respect to
$\arg z\in[-\pi,\pi]$. The second presentation in our main result \eqref{Basympcombined} has
the advantage that it is obvious that there are no issues when $a$ approaches a positive integer.
\end{proof}

In Table \ref{tab:combined}, we do illustrate that the implied constant in the order
estimate \eqref{Basympcombined} seems very reasonable in the full $z$-sector.
Note that the final result in Table \ref{tab:combined}, is for $z$ very close to
the boundary of the sector $|\arg(1+z)|\leq\pi$.

\begin{table}[ht]\renewcommand{\arraystretch}{1.2}
\begin{center}
\begin{tabular}{|c|c|c|}
\hline
 & $z=5$, $M_{\rm opt}=9$
 & $z=10$, $M_{\rm opt}=13$ \\
\hline
$\Phi(-z,s,a)$
& $1.3421782$
& $1.0889334$\\

approx \eqref{Basympcombined}
& $1.3421692$
& $1.0889332$\\

$\bigl|z^{N+1}R_N(z)\bigr|$
& $0.140$
& $0.158$\\
\hline
 & $z=100$, $M_{\rm opt}=27$
 & $z=1000$, $M_{\rm opt}=40$ \\
\hline
$\Phi(-z,s,a)$
& $0.50810464209509$
& $0.2350035297496389971$\\

approx \eqref{Basympcombined}
& $0.50810464209489$
& $0.2350035297496389969$\\

$\bigl|z^{N+1}R_N(z)\bigr|$
& $0.200$
& $0.197$\\
\hline
 & $z=10\iunit$, $M_{\rm opt}=16$
 & $z=-10+0.01\iunit$, $M_{\rm opt}=22$ \\
\hline
$\Phi(-z,s,a)$
& $0.98125249 - 0.54864116\iunit$
& $0.52526675 - 1.04285831\iunit$\\

approx \eqref{Basympcombined}
& $0.98125270 - 0.54864133\iunit$
& $0.52526654 - 1.04285810\iunit$\\

$\bigl|z^{N+1}R_N(z)\bigr|$
& $0.269$
& $0.297$\\

\hline
\end{tabular}
\caption{Approximation \eqref{Basympcombined} for the case $a=0.3$, $s=\frac34$
and $N=5$.}
\label{tab:combined}
\end{center}
\end{table}%

\section{A factorial series expansion}\label{Sect:factorial}
In this section, we will assume that $\bigl|1-\expe^{-2\pi\iunit a}\bigr|>1$, which is the case when,
for example, $a\in\bigl(\frac16,\frac56\bigr)$. We present $g(t)$ as
\begin{gather*}
 g(t)=\frac{1}{\expe^{\pi\iunit(a-t)}-\expe^{\pi\iunit(t-a)}}
 =\frac{\expe^{\pi\iunit(a-t)}}{\bigl(\expe^{2\pi\iunit a}-1\bigr)
 \bigl(1-\frac{1-\expe^{-2\pi\iunit t}}
 {1-\expe^{-2\pi\iunit a}}\bigr)}.
\end{gather*}
With the above assumption on $a$ we guarantee that the geometric progression
\begin{equation}\label{g3}
 g(t)=\expe^{\pi\iunit(a-t)}\sum_{n=0}^\infty
 \frac{\expe^{2\pi\iunit an}}{\bigl(\expe^{2\pi\iunit a}-1\bigr)^{n+1}}
 \bigl(1-\expe^{-2\pi\iunit t}\bigr)^n,
\end{equation}
converges uniformly for $t$ along the contour in \eqref{Bfact1}.
Using \eqref{g3} in \eqref{Bcontr2} gives us the expansion
\begin{equation}\label{Bfact1}
 B(z,s,a)=\frac{\expe^{\pi\iunit a}}{z^{a}}\sum_{n=0}^\infty
 \frac{\expe^{2\pi\iunit an}}{\bigl(\expe^{2\pi\iunit a}-1\bigr)^{n+1}}
 \int_{-\iunit\infty}^{(0+)}\expe^{(\ln(z)-\pi\iunit)t}t^{-s}\bigl(1-\expe^{-2\pi\iunit t}\bigr)^n\id t,
\end{equation}
which after taking $\tau=2\pi\iunit t$ and $x=\frac12-\frac{\ln z}{2\pi\iunit}$ becomes
\begin{equation}\label{Bfact2}
 B(z,s,a)=\frac{\expe^{\pi\iunit a}(2\pi\iunit)^{s-1}}{z^{a}}
 \sum_{n=0}^\infty
 \frac{\expe^{2\pi\iunit an}}{\bigl(\expe^{2\pi\iunit a}-1\bigr)^{n+1}}
 p_n(x,s),
\end{equation}
with
\begin{equation}\label{pn1}
 p_n(x,s)=\int_{\infty}^{(0+)}\expe^{-x\tau}\tau^{-s}\bigl(1-\expe^{-\tau}\bigr)^n\id\tau.
\end{equation}
If we assume that $\operatorname{Re} s<1$, we can collapse the loop contour
\begin{equation*}
 p_n(x,s)=\bigl(\expe^{-2\pi\iunit s}-1\bigr)
 \int_0^{\infty}\expe^{-x\tau}\tau^{-s}\bigl(1-\expe^{-\tau}\bigr)^n\id\tau.
\end{equation*}
Convergent expansion \eqref{Bfact2} is the main result of this section.
It is a factorial-type expansion because $p_n(x,s)$
is a generalisation (a fractional integral) of
\begin{equation*}
 \int_0^{\infty}\expe^{-x\tau}\bigl(1-\expe^{-\tau}\bigr)^n\id\tau=B(n+1,x)=
 \frac{n!}{x(x+1)(x+2)\cdots(x+n)}.
\end{equation*}
We still have to discuss the large $x$ asymptotic behaviour of $p_n(x,s)$,
and an alternative method to evaluate this function.

We use the expansion
\[\tau^{-s}=\sum_{m=0}^\infty c_m\bigl(1-\expe^{-\tau}\bigr)^{m-s},\]
which is the generating function for its coefficients
$c_m$,\footnote{We can identify $c_m$ in terms of generalised
Bernoulli coefficients: $c_m =(-1)^{m+1} s\frac{B_m^{(m-s)}}{(m - s)m!}$.
See \cite[\href{https://dlmf.nist.gov/24.16.E4}{equation~(24.16.4)}]{NIST:DLMF}).} and obtain
\begin{align}
 p_n(x,s)&=\bigl(\expe^{-2\pi\iunit s}-1\bigr)\sum_{m=0}^\infty c_m
 \int_0^{\infty}\expe^{-x\tau}\bigl(1-\expe^{-\tau}\bigr)^{n+m-s}\id\tau\nonumber\\
 &=\bigl(\expe^{-2\pi\iunit s}-1\bigr)\sum_{m=0}^\infty c_m B(n+m-s+1,x)\nonumber\\
 &\sim \bigl(\expe^{-2\pi\iunit s}-1\bigr)\Gamma(n-s+1)x^{s-n-1},\label{pn4}
\end{align}
as $x\to\infty$. Hence, $p_0(x,s),p_1(x,s),p_2(x,s),\ldots$ is definitely
an asymptotic sequence. These infinite series are conditionally convergent!

The binomial expansion of $\bigl(1-\expe^{-\tau}\bigr)^n$ in \eqref{pn1} will give us
\begin{equation}\label{pn5}
 p_n(x,s)=\frac{-2\pi\iunit\expe^{-\pi\iunit s}}{\Gamma(s)}
 \sum_{m=0}^n {n\choose m}(-1)^m(x+m)^{s-1},
\end{equation}
a finite sum. For modestly large $z$, our $x$ will not be large at all, and there is no
issue computing the $p_n(x,s)$ via~\eqref{pn5}. In this case, \eqref{Bfact2} will still converge,
but it can not be regarded an asymptotic expansion. However, in the case that $x$ is large
and comparing the terms in~\eqref{pn5} with estimate \eqref{pn4} it follows that the terms
in \eqref{pn5} are of the wrong size, hence cancellations. We can
deal with these cancellations via the identity
\[(1+\delta)^{\ell}=
\sum_{k=0}^{n-1}{\ell\choose k}\delta^k+\frac{(-\ell)_n(-\delta)^n}{n!}
\genhyperF{2}{1}{n-\ell}{1}{n+1}{-\delta}.\]
 This leads to
\begin{equation*}
 p_n(x,s)=\frac{-2\pi\iunit\expe^{-\pi\iunit s}}{\Gamma(s-n)}x^{s-n-1}
 \sum_{m=0}^n \frac{(-1)^m m^n}{m! (n-m)!}\genhyperF{2}{1}{n-s+1}{1}{n+1}{-m/x}.
\end{equation*}
Now the terms are of the correct size.

\subsection*{Acknowledgements}
This research was supported by a research Grant 60NANB20D126 from the National Institute
of Standards and Technology. The author thanks the referees for very helpful comments and
suggestions for improving the presentation.

\pdfbookmark[1]{References}{ref}
\LastPageEnding


\begin{thebibliography}{10}
\footnotesize\itemsep=0pt

\bibitem{BHNOD2018}
Bennett T., Howls C.J., Nemes G., Olde~Daalhuis A.B., Globally exact
 asymptotics for integrals with arbitrary order saddles, \href{https://doi.org/10.1137/17M1154217}{\textit{SIAM~J.~Math.
 Anal.}} \textbf{50} (2018), 2144--2177, \href{https://arxiv.org/abs/1710.10073}{arXiv:1710.10073}.

\bibitem{Cai2019}
Cai X.S., L\'opez J.L., A note on the asymptotic expansion of the {L}erch's
 transcendent, \href{https://doi.org/10.1080/10652469.2019.1627530}{\textit{Integral Transforms Spec. Funct.}} \textbf{30} (2019),
 844--855, \href{https://arxiv.org/abs/1806.01122}{arXiv:1806.01122}.

\bibitem{Erdelyi:1953:HTF1}
Erd\'elyi A., Magnus W., Oberhettinger F., Tricomi F.G., Higher transcendental
 functions. {V}ol.~{I}, Robert E.~Krieger Publishing Co., Inc., Melbourne, FL,
 1981.

\bibitem{Ferreira2004}
Ferreira C., L\'opez J.L., Asymptotic expansions of the {H}urwitz--{L}erch zeta
 function, \href{https://doi.org/10.1016/j.jmaa.2004.05.040}{\textit{J.~Math. Anal. Appl.}} \textbf{298} (2004), 210--224.

\bibitem{Navas2013}
Navas L.M., Ruiz F.J., Varona J.L., Asymptotic behavior of the {L}erch
 transcendent function, \href{https://doi.org/10.1016/j.jat.2012.03.006}{\textit{J.~Approx. Theory}} \textbf{170} (2013),
 21--31.

\bibitem{Nemes2017}
Nemes G., Error bounds for the asymptotic expansion of the {H}urwitz zeta
 function, \href{https://doi.org/10.1098/rspa.2017.0363}{\textit{Proc.~A.}} \textbf{473} (2017), 20170363, 16~pages,
 \href{https://arxiv.org/abs/1702.05316}{arXiv:1702.05316}.

\bibitem{NCA2023}
Nicholson M.D., Cheek D., Antal T., Sequential mutations in exponentially
 growing populations, \href{https://doi.org/10.1371/journal.pcbi.1011289}{\textit{PLOS Comput. Biol.}} \textbf{19} (2023),
 e1011289, 32~pages, \href{https://arxiv.org/abs/2208.02088}{arXiv:2208.02088}.

\bibitem{OD98}
Olde~Daalhuis A.B., Hyperasymptotic solutions of higher order linear
 differential equations with a singularity of rank one, \href{https://doi.org/10.1098/rspa.1998.0145}{\textit{R.~Soc. Lond.
 Proc. Ser.~A Math. Phys. Eng. Sci.}} \textbf{454} (1998), 1--29.

\bibitem{Olver1997}
Olver F.W.J., Asymptotics and special functions, \textit{AKP Class.}, A~K~Peters, Ltd.,
 \href{https://doi.org/10.1201/9781439864548}{Wellesley}, MA, 1997.

\bibitem{NIST:DLMF}
Olver F.W.J., Olde~Daalhuis A.B., Lozier D.W., Schneider B.I., Boisvert R.F.,
 Clark C.W., Miller B.R., Saunders B.V., Cohl H.S., McClain M.A., NIST digital
 library of mathematical functions, {R}elease 1.1.12 of 2023-12-15, available
 at \url{https://dlmf.nist.gov/}.

\bibitem{paris2016stokes}
Paris R.B., The {S}tokes phenomenon and the {L}erch zeta function,
 \textit{Math. Aeterna} \textbf{6} (2016), 165--179, \href{https://arxiv.org/abs/1602.00099}{arXiv:1602.00099}.

\bibitem{Paris2001}
Paris R.B., Kaminski D., Asymptotics and {M}ellin--{B}arnes integrals,
 \textit{Encyclopedia Math. Appl.}, Vol.~85, \href{https://doi.org/10.1017/CBO9780511546662}{Cambridge University Press},
 Cambridge, 2001.

\bibitem{Srivastava2011}
Srivastava H.M., Saxena R.K., Pog\'any T.K., Saxena R., Integral and
 computational representations of the extended {H}urwitz--{L}erch zeta
 function, \href{https://doi.org/10.1080/10652469.2010.530128}{\textit{Integral Transforms Spec. Funct.}} \textbf{22} (2011),
 487--506.

\end{thebibliography}
\end{document}